\documentclass[12pt]{article}

\usepackage{amsmath,amsfonts,amssymb}

\newenvironment{demo}{\noindent {\sl Proof}. \ }{\qed \bigskip}

\newtheorem{stheorem}{Theorem }[subsection]

\newtheorem{sprop}[stheorem]{Proposition }
\newtheorem{slemma}[stheorem]{Lemma }

\numberwithin{equation}{section}

\def\half{\frac{1}{2}}
\def\noi{\noindent}
\def\R{\mathbb{R}}
\def\N{\mathbb{N}}
\def\ep{\varepsilon}

\def\omb{\overline \omega}
\def\al{\alpha}

\def\de{\delta}

\def\tix{\tilde x}
\def\esp{{\rm I\mskip-4mu E}}
\def\Hess{{\rm Hess \ }}

\def\qed{\hbox{$\vcenter{\vbox{
   \hrule height 0.4pt\hbox{\vrule width 0.4pt height 6pt
    \kern5pt\vrule width 0.4pt}\hrule height 0.4pt}}$}}

\begin{document}

\title{
\bigskip \bigskip
  The monotonicity condition for Backward Stochastic Differential
  Equations on Manifolds}
\author{
\begin{tabular}{c}
Fabrice Blache \\
{\it{Institute for Applied Mathematics}}                 \\
{\it{University of Bonn, Poppelsdorfer Allee 82, 53115 Bonn}}  \\
{\it{Germany}}  \\
{\tt{E-mail: blache@wiener.iam.uni-bonn.de}}              \\
\end{tabular}}
\date{October 2005}

\maketitle

\vskip .5cm

\centerline{\small{\bf Abstract}}

\begin{center}
\begin{minipage}[c]{330pt}
{\small In \cite{blache03} and \cite{blache05}, we studied the problem of the existence and uniqueness of a solution to some general BSDE on manifolds. In these two articles, we assumed some Lipschitz conditions on the drift $f(b,x,z)$. The purpose of this article is to extend the existence and uniqueness results under weaker assumptions, in particular a monotonicity condition in the variable $x$. This extends well-known results for Euclidean BSDE.}
\end{minipage}
\end{center}

\bigskip
%%%%%%%%%%%%%%%%%%%%%%%%%%%%%%%%%%%%%%%%%%%%%%%%%%%%%%%%%%%%%%%%%%
%%%%%%%%%%%%%%%%%%%%%%   INTRODUCTION   %%%%%%%%%%%%%%%%%%%%%%%%%%
%%%%%%%%%%%%%%%%%%%%%%%%%%%%%%%%%%%%%%%%%%%%%%%%%%%%%%%%%%%%%%%%%%

\section{Reminder of the problem}
\label{par1}
Unless otherwise stated, we shall work on a fixed finite time interval
$[0;T]$; moreover, $(W_t)_{0 \le t \le T}$ will always denote a Brownian
Motion (BM for short) in $\R^{d_w}$, for a positive integer $d_w$. Moreover,
Einstein's summation convention will be used for repeated indices in
lower and upper position.

Let $(B_t^y)_{0 \le t \le T}$ denote the $\R^d$-valued
diffusion which is the unique strong solution of the following SDE :
\begin{equation}
\left\{
\begin{array}{rcl}
dB_t^y & = & b(B_t^y) dt + \sigma(B_t^y) d W_t  \\
 B_0^y & = & y,
\end{array}
\right.
\label{sde2}
\end{equation}
where $\sigma~: \R^d \rightarrow \R^{d \times d_w}$ and
$b : \R^d \rightarrow \R^d$ are $C^3$ bounded functions with bounded
partial derivatives of order 1, 2 and 3.

Let us recall the problem studied in \cite{blache03} and \cite{blache05}.
We consider a manifold $M$ endowed with a connection $\Gamma$, which
defines an exponential mapping. On $M$, we study the uniqueness and
existence of a solution to the equation (under infinitesimal form)
$$
(M+D)_0
\left\{
      \begin{array}{l}
           X_{t+ d t}= {\rm exp}_{X_t}(Z_t d W_t + f(B_t^y, X_t, Z_t) dt)\\
           X_T=U                         \\
      \end{array}
      \right.
$$
where
$Z_t \in {\cal L}(\R^{d_w}, T_{X_t}M)$ and $f(B_t^y,X_t,Z_t) \in T_{X_t}M$.

For details about links with PDEs, the reader is referred to the introductions of \cite{blache03} and \cite{blache05}.

\bigskip
In local coordinates $(x^i)$, the equation $(M+D)_0$ becomes the following
backward stochastic differential equation (BSDE in short)
$$(M+D)
\left\{
      \begin{array}{l}
          d X_t = Z_t d W_t + \left( - \half \Gamma_{jk}(X_t)
            ([Z_t]^k \vert [Z_t]^j) + f(B_t^y,X_t,Z_t) \right) d t  \\
      X_T=U.                                               \\
      \end{array}
    \right.   $$
We keep the same notations as in \cite{blache03}~: $( \cdot \vert \cdot)$
is the usual inner product in an Euclidean space, the summation convention
is used, and $[A]^i$ denotes the $i^{th}$ row of any matrix $A$; moreover,
\begin{equation}
\Gamma_{jk} (x) =
\left(
   \begin{array}{c}
       \Gamma^1_{jk}(x)    \\
       \vdots              \\
       \Gamma^n_{jk}(x)
   \end{array}
\right)
\label{christof2}
\end{equation}
is a vector in $\R^n$, whose components are the Christoffel symbols
of the connection. We keep the notations $Z_t$ for a matrix in
$\R^{n \times d_w}$ and $f$ for a mapping from
$\R^d \times \R^n \times \R^{n \times d_w}$ to $\R^n$.
The process $X$ will take its values in a compact set, and a solution of
equation $(M+D)_0$ will be a pair of processes
$(X,Z)$ in $M \times (\R^{d_w} \otimes TM)$ such that $X$ is continuous and
$\esp (\int_0^T \Vert Z_t \Vert_r^2 dt) < \infty$ for a
Riemannian norm $\Vert \cdot \Vert_r$;
in global coordinates $O \subset \R^n$,
$(X,Z) \in O \times \R^{n \times d_w}$ and
$\esp (\int_0^T \Vert Z_t \Vert ^2 dt) < \infty$ (see below for the
definitions of the norms).

We gave in \cite{blache03} and \cite{blache05} existence and
uniqueness results for the solutions of the BSDE $(M+D)$ for a drift
verifying some geometrical Lipschitz condition in the variables
$(b,x,z)$. Here we propose to weaken these assumptions, replacing
the Lipschitz condition in the variable $x$ by a monotonicity
condition and some other mild conditions; see Section \ref{par3}
below. The results for this condition were announced in
\cite{blache05}.

In Section \ref{par2}, we recall the general hypothesis. We give the
monotonicity condition and other assumptions on the drift $f$ in
Section \ref{par3}, as well as existence and uniqueness theorems
under the current assumptions. In Section \ref{par4}, we sketch the
changes to make in the uniqueness proofs in \cite{blache03} and
\cite{blache05}, to get uniqueness in the new context, and in
Section \ref{par5}, we prove the existence results. At the end in
Section \ref{par6}, we give the corresponding results for random
terminal times.

%%%%%%%%%%%%%%%%%%%%%%%%  general hypothesis  %%%%%%%%%%%%%%%%%%%%%%%%%%%
%%%%%%%%%%%%%%%%%%%%%%%%%%%%%%%%%%%%%%%%%%%%%%%%%%%%%%%%%%%%%%%%%%%%%%%%%
\section{Notations and hypothesis}
\label{par2}
In all the article, we suppose that a filtered probability space
$(\Omega, {\cal F}, P, ({\cal F}_t)_{0 \le t \le T})$ (verifying the usual
conditions) is given (with $T<\infty$ a deterministic time) on which
$(W_t)_t$ denotes a $d_w$-dimensional BM. Moreover, we
always deal with a complete Riemannian manifold  $M$ of dimension $n$,
endowed with a linear symmetric (i.e. torsion-free) connection whose
Christoffel symbols $\Gamma^i_{jk}$ are smooth; the connection does not
depend {\it a priori} on the Riemannian structure.

On $M$, $\delta$ denotes the Riemannian distance;
$\vert u \vert_r$ is the Riemannian norm for a tangent vector $u$
and $\vert u' \vert$ the Euclidean norm for a vector $u'$ in $\R^n$.
If $h$ is a smooth real function defined on $M$, its differential is
denoted by $D h$ or $h'$; the Hessian $\Hess h(x)$ is a bilinear
form the value of which is denoted by $\Hess h(x)<u,\overline u>$, for
tangent vectors (at $x$) $u$ and $\overline u$.\\
For $\beta \in \N^*$, we say that a function is $C^\beta$ on a closed set
$F$ if it is $C^\beta$ on an open set containing $F$.
Recall also that a real function $\chi$ defined on $M$ is said to be
convex if for any $M$-valued geodesic $\gamma$, $\chi \circ \gamma$ is
convex in the usual sense (if $\chi$ is smooth, this is equivalent to
require that $\Hess \chi$ be nonnegative).\\
For a matrix $z$ with $n$ rows and $d$ columns, $^t z$ denotes
its transpose,
$$\Vert z \Vert
= \sqrt{ \hbox{Tr} (z {}^t z)}
= \sqrt{\sum_{i=1}^d \vert [{}^t z]^i \vert ^2}$$
(${\rm Tr}$ is the trace of a square matrix) and
$\Vert z \Vert_r = \sqrt{\sum_{i=1}^d \vert [{}^t z]^i \vert_r ^2}$
where the columns of $z$ are considered as tangent vectors.
The notation $\Psi(x,x') \approx \delta(x,x')^\nu$ means that there is
a constant $c>0$ such that
$$ \forall x,x', \
\frac{1}{c} \ \delta(x,x')^\nu \le \Psi(x,x') \le c \ \delta(x,x')^\nu.$$

Throughout the article, we
consider an open set $O$ of $M$, relatively
compact in a local chart and an open set $\omega \not= \emptyset$
relatively compact in $O$, verifying that
\medskip
\\ $\bullet$ There is a unique geodesic in $\overline O$, linking any two
              points of $\overline O$, and depending smoothly on its
          endpoints;
\\ $\bullet$ $\omb = \{ \chi \le c \}$, the sublevel set of a smooth
              convex function $\chi$ defined on  $O$.
\medskip
\noi Note that $O$ will be as well considered as a subset of $\R^n$.

\noi We recall that in the case of a general connection, any point
$x$ of $M$ has a neighborhood $O$ for which the first property
holds; and when the Levi-Civita connection is used, the first
property is also true for regular geodesic balls.

%%%%%%%%%%%%%%%%%%%%%%%%%%%%%%%%%%%%%%%%%%%%%%%%%%%%%%%%%%%%%%%%%%%%%%%%%
%%%%%%%%%%%%%%%%%%%%%%%%%%%%%%%%%%%%%%%%%%%%%%%%%%%%%%%%%%%%%%%%%%%%%%%%%
%%%%%%%%%%%%%%%%%%%%%%%%%%%%%%%%%%%%%%%%%%%%%%%%%%%%%%%%%%%%%%%%%%%%%%%%%

\section{Case of a drift $f$ verifying a monotonicity condition}
\label{par3}

\subsection{Assumptions on $f$}
In this section, we explicit the assumptions on the drift $f$ under which
we will work; the main difference with \cite{blache03} and \cite{blache05}
is that we replace the Lipschitz property in $x$ by a monotonicity condition.
First, we suppose that $f$ verifies the Lipschitz condition in the
variables $b$ and $z$
\begin{eqnarray}
&& \exists L>0 , \ \forall b,b' \in \R^d,
  \forall x \in O,
  \forall z,z' \in {\cal L}(\R^{d_w}, T_{x}M),   \nonumber \\
&& \left\vert f(b,x,z) - f(b',x,z') \right\vert_r
\le L \Bigg( \vert b-b' \vert (1 + \Vert z \Vert_r + \Vert z' \Vert_r)
   + \left\Vert z - z' \right\Vert_r   \Bigg).
\label{lipz}
\end{eqnarray}
Before introducing the monotonicity condition in $x$, we recall
that, in the case of a general connection, $\Psi$ is a smooth,
convex and nonnegative function, vanishing on the diagonal only and
such that $\Psi \approx \delta^p$. In the case of the Levi-Civita
connection, we take $\Psi=\delta^2/2$
(in fact it is equivalent to take $\Psi = \delta$ off the diagonal).\\
The monotonicity condition is written
\begin{eqnarray}
&& \exists \nu \in \R , \ \forall b \in \R^d,
  \forall x, x' \in O,
  \forall z \in {\cal L}(\R^{d_w}, T_{x}M), \nonumber \\
&& D \Psi (x,x') \cdot
     \left(
        \begin{array}{c}
          f(b,x,z)  \\
          f(b,x', \overset{x'}{\underset{x}{\Vert}} z)
        \end{array}
     \right)
 \ge
 \nu  \Psi(x,x') (1 + \Vert z \Vert_r)
 \label{monot}
\end{eqnarray}
where
$$\overset{x'}{\underset{x}{\Vert}}  z$$ denotes
the parallel transport (defined by the connection) of the $d_w$ columns
of the matrix $z$ (considered as tangent vectors) along the unique
geodesic between $x$ and $x'$.
This monotonicity condition replaces here the well-known
monotonicity condition involving the inner product in an Euclidean
space (see e.g. Assumption $(4)$ in \cite{darlpard97} or Assumption
$(H3)$ in \cite{briandhupard03}).
Note that here we have a lower bound on $D \Psi$
(and not an upper bound as in the articles cited above) because in the
equation $(M+D)$, the drift $f$ is given with a ``plus" sign.

We also need the following uniform boundedness condition
\begin{equation}
\exists L_2>0 , \forall b \in \R^d, \forall x \in O, \vert f(b,x,0)
\vert_r \le L_2
\label{upperboundf}
\end{equation}
and the continuity in the $x$ variable~:
\begin{equation}
\forall b \in \R^d, \ \forall z \in {\cal L}(\R^{d_w}, T_{x}M),
\ \ x \mapsto f(b,x,z) \ \hbox {is continuous}.
\label{continux}
\end{equation}

\subsection{The results}
We state now the existence and uniqueness theorem, which generalizes
Theorem 1.4.1 of \cite{blache03} and Theorem 1.3.1 of
\cite{blache05} to our context. We need the following assumption~:
$$(H) \ \ f \hbox{ is pointing outward on the boundary of } \omb.$$

\begin{stheorem}
\label{existence}
We consider the BSDE $(M+D)$ with terminal random variable
$U \in \omb= \{ \chi \le c \}$. We suppose
that $f$ verifies conditions \eqref{lipz}, \eqref{monot},
\eqref{upperboundf} and \eqref{continux},
and that $\chi$ is strictly convex (i.e. $\Hess \chi$ is
positive definite). Moreover for $(i)$ and $(ii)$, $f$ is also
supposed to verify $(H)$.
\begin{item}
(i) If f does not depend on $z$, the BSDE has a unique solution
$(X_t,Z_t)_{0 \le t \le T}$ such that $X$ remains in $\omb$.
\end{item}
\begin{item}
(ii) If the Levi-Civita connection is used, then the BSDE has yet a unique solution
$(X_t,Z_t)_{0 \le t \le T}$ with $X$ in $\omb$ (this is true in particular for regular geodesic balls) .
\end{item}
\begin{item}
(iii) In the case of a general connection, each point $q$ of $M$ has
a neighborhood $O_q \subset O$ such that, if $\omb \subset O_q$ and
$f$ verifies hypothesis $(H)$, then the BSDE $(M+D)$ has a unique
solution $(X_t,Z_t)_{0 \le t \le T}$ such that $X$ remains in
$\omb$. The neighborhood $O_q$ depends on the geometry of the
manifold, but not on the drift $f$.
\end{item}
\end{stheorem}

\subsection{Preliminary lemmas}
First we recall two results~: Proposition 2.2.1 in \cite{blache05}
and Lemma 3.2.1 of \cite{blache03} (in fact Lemma \ref{upbddpsi}
below is a straightforward consequence of this lemma and Proposition
\ref{proptp}).

\begin{sprop}
\label{proptp}
Here $O$ is considered, via a system of local coordinates, as an open
subset of $\R^n$. There is a $C>0$ such that for every
$(x,x') \in O \times O$ and  $(z,z') \in T_xM \times T_{x'}M$, we have
$$
\left\vert \overset{x'}{\underset{x}{\Vert}} z - z' \right\vert_r
 \le C
  \left( \vert z-z' \vert + \delta(x,x')(\vert z \vert + \vert z' \vert)
       \right)
$$
and
\begin{equation}
  \vert z-z' \vert \le C
\left( \left\vert \overset{x'}{\underset{x}{\Vert}}
                                    z - z' \right\vert_r
  + \delta(x,x')(\vert z \vert_r + \vert z' \vert_r)
       \right).
 \label{2tp2}
\end{equation}
\end{sprop}

\begin{slemma}
\label{upbddpsi} Suppose that $\Psi(x,x') \approx \delta(x,x')^p$ on
$O \times O$ where $p$ is an even positive integer (since $\Psi$ is
smooth). Then there is $C>0$ such that, for all vectors
$(z,z') \in
T_x M \times T_{x'} M$,
\begin{equation}
 \left\vert D\Psi (x,x') \cdot
                                   \left(
                                 \begin{array}{c}
                                     z \\
                                     z'
                                  \end{array}
                                   \right)
 \right\vert \le
 C  \delta(x,x')^{p-1} \left( \delta(x,x') (\vert z \vert_r
                                + \vert z' \vert_r)
       + \vert \overset{x'}{\underset{x}{\Vert}} z - z'\vert_r \right).
\label{majdpsi}
\end{equation}
\end{slemma}

Then we have the following lemma, which is in fact the main argument in
order to prove uniqueness. For notational convenience, we denote $f(b,x,z)$ by $f$ and $f(b,x',z')$ by $f'$.
\begin{slemma}
\label{boundDpsi} There is $\hat C>0$ such that, for all $b$ in
$\R^d$, $x,x'$ in $O$, $z \in {\cal L}(\R^{d_w}, T_{x}M)$ and $z'
\in {\cal L}(\R^{d_w}, T_{x'}M)$,
\begin{equation}
 D\Psi (x,x') \cdot
                                   \left(
                                 \begin{array}{c}
                                     f \\
                                     f'
                                  \end{array}
                                   \right)
 \ge
 - \hat C  \delta(x,x')^{p-1}
 \left(\left\Vert \overset{x'}{\underset{x}{\Vert}} z - z' \right\Vert_r          + \delta (x,x') (1 + \Vert z \Vert_r + \Vert z' \Vert_r) \right).
\end{equation}
\end{slemma}

\begin{demo}
Note that we use the same $b$.
We have
\begin{eqnarray*}
D\Psi (x,x') \cdot
                                   \left(
                                 \begin{array}{c}
                                     f \\
                                     f'
                                  \end{array}
                                   \right)
& = &
D\Psi (x,x') \cdot
                      \left(
                        \begin{array}{c}
                          f(b,x,z) \\
                          f(b,x',\overset{x'}{\underset{x}{\Vert}} z)
                        \end{array}
                      \right) \\
&   &
+ D\Psi (x,x') \cdot
             \left(
               \begin{array}{c}
                0 \\
                f(b,x',z')-f(b,x',\overset{x'}{\underset{x}{\Vert}} z)
              \end{array}
             \right).\\
&\ge&
\nu \Psi(x,x') (1+ \Vert z \Vert) \\
&   & - C  \delta(x,x')^{p-1}
 \left(\left\Vert \overset{x'}{\underset{x}{\Vert}} z - z' \right\Vert_r
     + \delta(x,x') (1 + \Vert z \Vert_r + \Vert z' \Vert_r) \right).
\end{eqnarray*}
In the last inequality, we have used the monotonicity assumption
\eqref{monot} for the first term, and for the second term the
Lipschitz assumption \eqref{lipz} and Lemma \ref{upbddpsi} above.
This gives the result since $\Psi \approx \delta^p$.
\end{demo}

To end this part, we give the following result, which means that $f$
has at most a linear growth in the variable $z$. It will be useful
in the sequel.
\begin{slemma}
Under the above assumptions \eqref{lipz} and \eqref{upperboundf},
there is a constant $C$ such that
$$
\forall \ b,x,z, \ \
\vert f(b,x,z) \vert_r \le C ( \Vert z \Vert_r + 1).
$$
\label{lingro}
\end{slemma}

\begin{demo}
It is a straightforward consequence of the two assumptions
\eqref{lipz} and \eqref{upperboundf}~:
\begin{eqnarray*}
\forall \ b,x,z, \ \ \vert f(b,x,z) \vert_r
&\le& \vert f(b,x,z) - f(b,x,0) \vert_r + \vert f(b,x,0) \vert_r \\
&\le& L \Vert z \Vert_r + L_2.
\end{eqnarray*}
This gives the result.
\end{demo}

%%%%%%%%%%%%%%%%%%%%%%%%%%%%%%%%%%%%%%%%%%%%%%%%%%%%%%%%%%%%%%%%%%%%%%%%%
%%%%%%%%%%%%%%%%%%%%%%%%%%%%%%%%%%%%%%%%%%%%%%%%%%%%%%%%%%%%%%%%%%%%%%%%%
\section{The uniqueness property}
\label{par4} We recall first the method used in \cite{blache03} and
in \cite{blache05},  then we will sketch the main changes to make
under our new assumptions.

\subsection{General results}
Consider two solutions $(X_t,Z_t)_{0 \le t \le T}$ and
$(X'_t,Z'_t)_{0 \le t \le T}$  of $(M+D)$ such that $X$ and $X'$
remain in $\omb$ and $X_T=Y_T=U$ (we will say sometimes
``$\omb$-valued solutions of $(M+D)$"). Let
$$\tilde X_s = (X_s,X'_s) \ \ \hbox{ and }
   \tilde Z_s=  \left(
          \begin{array}{c}
             Z_s  \\
             Z'_s
          \end{array}  \right).$$
To prove uniqueness, the idea is to show that the process
$(S_t)_t=(\exp (A_t) \Psi (\tilde X_t))_t$
where
$$
A_t= \lambda t + \mu \int_0^t (\Vert Z_s \Vert_r^\al  +
                           \Vert Z'_s \Vert_r^\al ) ds,
$$
is a submartingale for appropriate nonnegative constants $\lambda$,
$\mu$ and $\al$, and a suitable function $\Psi$, smooth on $O \times
O$ (of course in general we will have to prove some integrability
properties).

If we apply It\^o's formula to $(S_t)$, we get

\begin{eqnarray}
e^{A_t}  \Psi (\tilde X_t) - \Psi (\tilde X_0)
& = &\int_0^t e^{A_s} d (\Psi (\tilde X_s))
+ \int_0^t e^{A_s} (\lambda + \mu (\Vert Z_s \Vert_r^\al + \Vert Z'_s \Vert_r^\al))
                       \Psi (\tilde X_s) ds               \nonumber \\
& = &\int_0^t e^{A_s} D \Psi (\tilde X_s)
          \left( \tilde Z_s d W_s \right)                 \nonumber \\
&   & + \half \int_0^t e^{A_s}
          \left(  \sum_{i=1}^{d_W} {}^t [{}^t \tilde Z_s]^i
   \Hess \Psi (\tilde X_s) [{}^t \tilde Z_s]^i \right) ds \nonumber \\
&  & + \int_0^t e^{A_s} D \Psi(\tilde X_s)
                                   \left(
                                 \begin{array}{c}
                                     f(B_s^y,X_s,Z_s)  \\
                                     f(B_s^y,X'_s,Z'_s)
                                  \end{array}
                                   \right) ds             \nonumber \\
&  & + \int_0^t e^{A_s} \Psi (\tilde X_s) (\lambda + \mu (\Vert Z_s \Vert_r^\al
     + \Vert Z'_s \Vert_r^\al)) ds.
\end{eqnarray}

Thus to prove the submartingale property, we need to show the nonnegativity of the sum
\begin{eqnarray}
\half \sum_{i=1}^{d_W} {}^t [{}^t \tilde Z_t]^i
        \Hess \Psi (\tilde X_t) [{}^t \tilde Z_t]^i
& + &  D \Psi(\tilde X_t)
                                   \left(
                                 \begin{array}{c}
                                     f(B_t^y,X_t,Z_t)  \\
                                     f(B_t^y,X'_t,Z'_t)
                                  \end{array}
                                   \right)             \nonumber     \\
& + & (\lambda + \mu (\Vert Z_t \Vert_r^\al + \Vert Z'_t \Vert_r^\al))
            \Psi(\tilde X_t).
\label{pos}
\end{eqnarray}

\subsection{Sketch of the proof}
{\it Case of a drift independent of $z$}

In this case, the result is straightforward; indeed, from Lemma \ref{boundDpsi}, there is $C>0$ such that,
for all $b$ in $\R^d$, $x,x'$ in $\omb$,

$$ D\Psi(x,x') \cdot
                                   \left(
                                 \begin{array}{c}
                                     f(b,x)  \\
                                     f(b,x')
                                  \end{array}
                                   \right)
 \ge -C \delta(x,x')^p \ge - \tilde C \Psi(x,x').
 $$
Therefore in order to make the sum \eqref{pos} nonnegative, it is
sufficient to take $\lambda$ larger than $\vert \tilde C \vert$ and
$\mu=0$.

\bigskip
\noi {\it Uniqueness for drifts depending on $z$}

As we remarked above, uniqueness is a consequence of the
nonnegativity  of the sum \eqref{pos}. Of course the estimates on
the Hessian term which are given in \cite{blache03} and
\cite{blache05} are still valid; thus we have only to verify that
the same kind of estimates hold for the term involving $f$ in
\eqref{pos}. Moreover, we need to keep the integrability conditions
which make integrable the process $(S_t)_t=(\exp (A_t) \Psi (\tilde
X_t))_t$.

 We begin with the integrability condition; in fact it is based upon
 Lemma 3.4.2 of \cite{blache03} which we recall now.
\begin{slemma}
\label{Emu}
Suppose that we are given a positive
constant $\al$ and a $C^2$ function $\phi$ on $\omb$ satisfying
$C_{min} \le \phi(x) \le C_{max}$ for some positive $C_{min}$ and
$C_{max}$. Suppose moreover that
$\Hess \phi + 2 \al \phi \le 0$ on $\omb$; this means that
\begin{equation}
  \Hess \phi(x)<u,u>+2 \al \phi(x) \vert u \vert_r^2 \le 0.
\end{equation}
Then, for every $\ep>0$, any $\omb$-valued solution of $(M+D)$ belongs
to $(\cal E_{\al - \ep})$.
\end{slemma}
An accurate examination of the proof given in \cite{blache03} shows
that the only hypothesis needed on the drift $f$ is the linear
growth in $z$,  and this is yet verified, according to Lemma
\ref{lingro}.

Therefore it remains to verify that the estimates on $D \Psi (x,x')(f,f')$
are still valid. We just check that this is  true in the three cases dealt
with in the two articles \cite{blache03} and \cite{blache05}.\\

In the case of a Levi-Civita connection, two different approaches
were developed~: \\
$(1)$ On the one hand for nonpositive sectional curvatures, the
function $\Psi$ is just $\de^2/2$ and Lemma \ref{boundDpsi} gives
the estimate
\begin{equation}
D\Psi (\tix)
                                   \left(
                                 \begin{array}{c}
                                     f(b,x,z)    \\
                                     f(b,x',z')
                                  \end{array}
                                   \right)
 \ge
 - C \delta^2(x,x')(1+\Vert z \Vert_r +
                                   \Vert z' \Vert_r)
       - \frac{1}{4} \left\Vert
          \overset{x'}{\underset{x}{\Vert}} z -z' \right\Vert_r^2.
\label{estim1}
\end{equation}
$(2)$ On the other hand, if positive sectional curvatures are allowed
(bounded by $K>0$), we defined in \cite{blache05}
$$\Psi(x,x')
  = \Psi_a(x,x')= \sin^a \left( \sqrt K \frac{\delta(x,x')}{2} \right)$$
and therefore, for $z \in {\cal L}(\R^{d_w}, T_xM)$ and
$z' \in {\cal L}(\R^{d_w}, T_{x'}M)$,
\begin{eqnarray}
D \Psi (\tix) \cdot     \left(
                                 \begin{array}{c}
                                     f(b,x,z) \\
                                     f(b,x',z')
                                  \end{array}
                                     \right)
& = & a \sin^{a-1}(y) \cos y \cdot \frac{\sqrt K}{2}
 \delta'(\tix) \cdot              \left(
                                 \begin{array}{c}
                                     f(b,x,z)  \\
                                     f(b,x',z')
                                  \end{array}
                                    \right)            \nonumber  \\
&\ge& - a \sin^{a-1}(y) \frac{\sqrt K}{2}              \nonumber\\
&   & \hskip 1cm \times L_0 \left( \delta(\tix)
            (1+\Vert z \Vert_r + \Vert z' \Vert_r )
      + \left\Vert \overset{x'}{\underset{x}{\Vert}} z - z'  \right\Vert_r
         \right)                                       \nonumber   \\
&\ge& - a \frac{\pi}{2} \Psi(\tix) L(1+ \Vert z \Vert_r + \Vert z'                \Vert_r )
      - C_1 \sin^{a-1}y  \left\Vert \overset{x'}{\underset{x}{\Vert}}
             z - z'  \right\Vert_r                       \nonumber\\
&\ge& - C \Psi(\tix) - \frac{e-1}{2} \frac{K}{4} \Psi(\tix)
           (\Vert z \Vert^2_r + \Vert z' \Vert^2_r)     \nonumber \\
&   & - \frac{\al}{2} \sin^{a-2} y
       \left\Vert \overset{x'}{\underset{x}{\Vert}} z - z'  \right\Vert_r^2.
\label{estim2}
\end{eqnarray}
for appropriate constants $a$, $e$ and $\alpha$ (which are in fact
determined according to the estimates upon $\Hess \Psi$ in the same
way as in \cite{blache05}), and positive constants $C_1$ and $C$.
Note that the first inequality is a consequence of the monotonicity
condition \eqref{monot} with $\Psi= \delta$ instead of $\delta^2$
(we have already remarked that we have equivalence off the
diagonal).

$(3)$ In the case of a more general connection, Lemma
\ref{boundDpsi}  gives, for any $\ep >0$,
\begin{eqnarray}
D \Psi \cdot
                                   \left(
                                 \begin{array}{c}
                                     f(b,x,z)  \\
                                     f(b,x',z')
                                  \end{array}
                                   \right)
& \ge &
 - C_\ep \Psi(x,x') (1 + \Vert z \Vert_r + \Vert z' \Vert_r)  \nonumber \\
&     & - \ep  \delta^{p-2}(x,x') \left\Vert
     \overset{x'}{\underset{x}{\Vert}} z - z' \right\Vert_r^2.
\label{estim3}
\end{eqnarray}

It turns out that the estimate \eqref{estim1} (resp. \eqref{estim2},
resp. \eqref{estim3}) is exactly what is needed on the term $D
\Psi(f,f')$ in order to get the nonnegativity of the sum
\eqref{pos}; thus it suffices  to replicate the proof of Lemma 3.4.5
in \cite{blache03} (resp. Proposition 3.3.2 in \cite{blache05},
resp. Theorem 3.2.2 in \cite{blache05}). This concludes the
uniqueness part.

%%%%%%%%%%%%%%%%%%%%%%%%%%%%%%%%%%%%%%%%%%%%%%%%%%%%%%%%%%%%%%%%%%%%%%%%%
%%%%%%%%%%%%%%%%%%%%%%%%%%%%%%%%%%%%%%%%%%%%%%%%%%%%%%%%%%%%%%%%%%%%%%%%%
%%%%%%%%%%%%%%%%%%%%%%%%%%%%%%%%%%%%%%%%%%%%%%%%%%%%%%%%%%%%%%%%%%%%%%%%%
\section{Existence property}
\label{par5}
First we give the framework of the proof.

Suppose that $c>0$ and $\chi$ reaches its infimum at $p \in \omega$
with $\chi(p)=0$. Then consider the following mapping, defined on a
normal open neighborhood of $\omb$ centered at $p$ (so it is in
particular a neighborhood of $0$ in $\R^n$)
$$ y \mapsto \frac{\sqrt{c_2} y}
                  {\sqrt{c_2 \Vert y \Vert^2 + c-\chi(y)} };$$
for $c_2>0$ small enough, it is a diffeomorphism from an open set
$O_1$ (relatively compact in $O$ and containing $\omb$) onto an open
neighborhood $N$ of $\overline {B(0,1)}$, such that $\omb$ is sent
onto $\overline {B(0,1)}$ (in fact, it is sufficient to take $c_2
\le \half \lambda_\chi$, where $\lambda_\chi$ denotes the (positive)
infimum on $O_1$ of the eingenvalues of $\Hess \chi$). Using this
diffeomorphism, we can work in a local chart $O$ (take $O:=N$) such
that $\omb = \overline {B(0,1)}$ and $\chi(0)=0$. We will mainly use
the Euclidean norm in this section; note in particular that
\eqref{lipz}, \eqref{monot} and \eqref{upperboundf} remain true if
the Riemannian norm is replaced by the Euclidean one
(changing accordingly the constants).\\
Remark that hypothesis $(H)$ means, in the new coordinates, that the
radial component $f^{\rm rad}(b,x,z)$ of $f(b,x,z)$ is nonnegative
for $x \in
\partial \omb$.

The goal in this part is to approximate $f$ by Lipschitz functions
for which the existence theorems of \cite{blache03} and of
\cite{blache05} apply. However, unlike in these articles, we cannot
use the convolution product directly; indeed, we don't have precise
estimates for $f$ in a neighborhood of $\omb$, therefore we can't
prove that these approximations be pointing outward on $\partial
\omb$. In order to give such estimates, we need first to introduce
new functions $f_k$ with which we shall use the convolution product.

\bigskip
{\it First step.} We truncate the growth of the drift $f$, both in
$b$ and $z$, in the following way. Let $\phi : \R \rightarrow \R_+$
be a smooth nonincreasing function with
$$
\begin{cases}
\phi (u)=1 & \ \hbox{ if } \  u \le 1 \\
\phi (u)=0 & \ \hbox{ if } \  u \ge 2; \\
\end{cases}$$
then put $\phi_k(u)= \phi(u-(k-1))$ for any positive integer $k$.
Note that $0 \le \phi_k \le 1$, with $\phi_k(u)=1$ if $u \le k$ and
$\phi_k(u)=0$ if $u \ge k+1$. Let
\begin{equation}
f_k(b,x,z)=f(b,x,z) \phi_k(\vert b \vert) \phi_k(\Vert z \Vert);
\label{deffk}
\end{equation}
in
particular, $f_k$ tends to $f$ when $k$ tends to infinity. We have
the

\begin{slemma}
\label{lipphi} There is a constant $C$, independent of $k$, such
that
\begin{eqnarray*}
(i)
& &
\forall x,x' \in O, \forall z \in \R^{n d_w}, \ \ \left\vert
\phi_k \left( \left\Vert \overset{x'}{\underset{x}{\Vert}} z
\right\Vert \right) - \phi_k (\Vert z \Vert) \right\vert \le C k
\delta(x,x');\\
(ii) & &
\forall x \in O, \forall z, z' \in \R^{n d_w}, \forall b
\in \R^d,\\
& &
\vert f(b,x,z') \vert \cdot \vert \phi_k (\Vert z' \Vert)-
\phi_k (\Vert z \Vert) \vert \le C (1+ k \wedge \Vert z' \Vert)
\Vert z-z' \Vert.
\end{eqnarray*}
\end{slemma}

\begin{demo}
The function $\phi_k$ verifies $ \vert \phi_k(u) - \phi_k(u') \vert
\le C_0 \vert u-u' \vert (1_{ \{ u \le k+1 \} }+ 1_{ \{ u' \le k+1\}
});$
using Proposition \ref{proptp}, we get
\begin{eqnarray*}
\left\vert \phi_k \left( \left\Vert
\overset{x'}{\underset{x}{\Vert}} z \right\Vert \right) - \phi_k
(\Vert z \Vert) \right\vert
&\le&
 C_0   \left\Vert \overset{x'}{\underset{x}{\Vert}} z - z\right\Vert
   (1_{ \{ \Vert z \Vert \le k+1\} } + 1_{ \{\Vert \overset{x'}
   {\underset{x}{\Vert}} z \Vert \le k+1\} })  \\
&\le&
C_1 \left( \delta (x,x') \Vert z \Vert 1_{ \{ \Vert z \Vert
\le k+1\} } + \delta (x,x') \Vert \overset{x'}{\underset{x}{\Vert}}
z \Vert
1_{ \{ \Vert \overset{x'}{\underset{x}{\Vert}} z \Vert \le k+1\} } \right) \\
&\le& C (k+1) \delta(x,x'),
\end{eqnarray*}
hence $(i)$.

For $(ii)$, as $\phi_k$ is Lipschitz (with a constant independent of
$k$), we have
$$\vert f(b,x,z') \vert \cdot \vert \phi_k
(\Vert z' \Vert)- \phi_k (\Vert z \Vert) \vert \le C (1+  \Vert z'
\Vert) \Vert z-z' \Vert.$$ As $f$ has a linear growth in $z$, it is
thus sufficient to prove
$$\Vert z' \Vert \cdot \vert \phi_k
(\Vert z' \Vert)- \phi_k (\Vert z \Vert) \vert \le C (1+ k) \Vert
z-z' \Vert.$$ The only case for which we have to be careful is $z'
\ge k+1$ and $z <k+1$. We have in this case
\begin{eqnarray*}
\Vert z' \Vert \cdot \vert \phi_k (\Vert z' \Vert)- \phi_k (\Vert z
\Vert) \vert
& \le&
(\Vert z'-z \Vert + \Vert z \Vert) \cdot
\vert \phi_k (\Vert z' \Vert)- \phi_k (\Vert z \Vert) \vert \\
&\le&
2 \Vert z'-z \Vert + (k+1) \Vert z'-z \Vert
\end{eqnarray*}
and that completes the proof.
\end{demo}

\begin{sprop}
\label{propfk}
We have the following properties for the new drift $f_k$~:
\begin{item}{(i)}
It verifies a Lipschitz property in $b$ and $z$ in the following sense~:
$$ \exists \hat L >0 , \ \forall b,b' \in \R^d,
  \forall x \in O,
  \forall z,z' \in \R^{n d_w}, $$
\begin{eqnarray*}
\left\vert f_k(b',x,z') - f_k(b,x,z) \right\vert &\le& \hat L \Bigg(
(\vert b-b' \vert + \Vert z - z' \Vert ) \times \\
&   &
  (1 + k \wedge \Vert z \Vert + k \wedge \Vert z' \Vert) \Bigg)
\end{eqnarray*}
\end{item}

\begin{item}{(ii)}
It verifies a monotonicity condition like \eqref{monot}, but now
with a constant depending on $k$~:
\begin{eqnarray*}
&& \exists C>0 , \ \forall b \in \R^d,
  \forall x, x' \in O,
  \forall z \in \R^{n d_w}, \nonumber \\
&& D \Psi (x,x') \cdot
     \left(
        \begin{array}{c}
          f_k(b,x,z)  \\
          f_k(b,x', \overset{x'}{\underset{x}{\Vert}} z)
        \end{array}
     \right)
 \ge
 -C k  \Psi(x,x') (1 + \Vert z \Vert).
\end{eqnarray*}
\end{item}

\begin{item}{(iii)}
When $z=0$, it is uniformly bounded as in \eqref{upperboundf} (independently of $k$).
\end{item}

\begin{item}{(iv)}
It is continuous in $x$ in the sense of \eqref{continux}.
\end{item}

\begin{item}{(v)}
It is pointing outward on the boundary of $\omb$.
\end{item}

\begin{item}{(vi)}
Recall from the beginning of this section that $\overline O_1$ is a
compact set verifying $\omb \subset O_1 \subset \overline O_1
\subset O$. Then $f_k$ is uniformly continuous in the $x$ variable,
i.e.
$$
\begin{array}{c}
\forall \ep >0, \exists \eta_k >0 \ \hbox{ s.t.} \ \forall (b,x,z) \in
   \R^d \times \overline O_1 \times \R^{n d_w}, \forall v \ \ s.t.
   \ \vert v \vert \le \eta_k,  \\
 \left\vert f_k(b,x,z)- f_k(b,x+v,z) \right\vert \le \ep.
 \end{array}
 $$
\end{item}
\end{sprop}

\begin{demo}
Using that $\phi_k(u) \le 1$, we have
\begin{eqnarray*}
\vert f_k(b',x,z') - f_k(b,x,z) \vert
&\le&
\vert f(b',x,z') \phi_k (\vert b' \vert) - f(b,x,z') \phi_k (\vert b \vert) \vert
 \cdot \phi_k(\Vert z' \Vert) \\
&   &
+ \vert f(b,x,z') \phi_k(\Vert z' \Vert) - f(b,x,z) \phi_k
(\Vert z \Vert)  \vert \\
&\le&
\vert f(b',x,z') \vert
\cdot \vert \phi_k (\vert b' \vert) -
\phi_k (\vert b \vert) \vert
\cdot \phi_k(\Vert z' \Vert) \\
&  &
+ \vert f(b',x,z') - f(b,x,z') \vert \cdot \phi_k(\Vert z'
\Vert) \\
&  &
+ \vert f(b,x,z') \vert \cdot
  \vert \phi_k(\Vert z'\Vert) -  \phi_k(\Vert z \Vert) \vert   \\
&  &
+ \vert f(b,x,z') - f(b,x,z) \vert \cdot \phi_k(\Vert z \Vert) \\
&\le&
C_1 (1+ \Vert z' \Vert) \vert b-b' \vert 1_{\Vert z' \Vert \le
k+1}\\
&  & + C_2 (1+ k \wedge \Vert z' \Vert) \Vert z-z' \Vert\\
&  & + C_3 \Vert z-z' \Vert.
\end{eqnarray*}
In the last inequality we have used \eqref{lipz}, Lemma \ref{lingro}
(for the Euclidean norm) and $(ii)$ of Lemma \ref{lipphi}. Since
$(1+ \Vert z' \Vert) 1_{\Vert z' \Vert \le k+1} \le C k \wedge \Vert
z' \Vert$, this proves $(i)$.

The properties $(iii)$, $(iv)$ and $(v)$ follow obviously from the
same properties for $f$. Moreover, $(vi)$ is just a consequence of
the continuity of $f$ on the compact set $\overline {B(0,k+1)}
\times \overline O_1 \times \overline {B(0,k+1)}$ (recall that if $b
\notin \overline {B(0,k+1)}$ or $z \notin \overline {B(0,k+1)}$,
then $f_k(b,x,z)=0$).

For $(ii)$, we write $f_k$ (resp $f$) for $f_k(b,x,z)$ (resp. $f(b,x,z))$ and $f_k'$ (resp $f'$) for
$f_k(b,x', \overset{x'}{\underset{x}{\Vert}} z)$
(resp. $f(b,x', \overset{x'}{\underset{x}{\Vert}} z)$)
and we have
\begin{eqnarray*}
D \Psi \cdot
     \left(
        \begin{array}{c}
          f_k \\
          f_k'
        \end{array}
     \right)
&=&
 D \Psi \cdot
     \left(
        \begin{array}{c}
          \phi_k(\Vert z \Vert) \phi_k(\vert b \vert) f\\
          \phi_k(\Vert \overset{x'}{\underset{x}{\Vert}} z \Vert) \phi_k(\vert b \vert) f'
        \end{array}
     \right) \\
&=&
D \Psi \cdot
     \left(
        \begin{array}{c}
          f \\
          f'
        \end{array}
     \right)
\phi_k(\vert b \vert) \phi_k(\Vert z \Vert)\\
& + &
 D \Psi \cdot
     \left(
        \begin{array}{c}
          0 \\
          (\phi_k(\Vert \overset{x'}{\underset{x}{\Vert}} z \Vert)-
          \phi_k(\Vert z \Vert)) f'
        \end{array}
     \right)
     \phi_k(\vert b \vert) \\
&\ge& (\nu \wedge 0) \Psi(x,x') (1+ \Vert z \Vert)
      - C_1 k \Psi(x,x') (1+ \Vert z \Vert) \\
&\ge& - C k \Psi(x,x') ( 1 + \Vert z \Vert).
\end{eqnarray*}
In the first inequality, we have used for the first term the
monotonicity assumption \eqref{monot} on $f$, and for the second
term the Euclidean form of Lemma \ref{upbddpsi}, $(i)$ in Lemma
\ref{lipphi} and the linear growth of $f$ in $z$ (Lemma
\ref{lingro}). Note that the constant $C>0$ at the end does not
depend on $k$. Then we have $(ii)$.
\end{demo}

The interest of the new drifts $f_k$ lies in condition $(vi)$;
indeed, this gives some information about the behavior of $f_k$
around the boundary $\partial \omb$. Then now we can construct new
approximations which will verify Lipschitz conditions in all the
variables.

\bigskip
{\it Second step.} Extend each mapping $f_k$ to $\R^d \times \R^n
\times \R^{n d_w}$ by putting $f_k(b,x,z)=0$ if $x \notin O$ and
define (on $\R^d \times \R^n \times \R^{n d_w}$) the convolution
product for $l \in \N^*$ $f_{k,l} = f_k * \rho_l$~:
$$ f_{k,l}(b,x,z)
    = \int_{\R^d \times \R^n \times \R^{n d_w}}
        f_k((b,x,z)-(\beta,y,w)) \rho_l(\beta,y,w) d (\beta,y,w)$$
where $\rho_l(b,x,z) = l \rho (l \Vert (b,x,z) \Vert)$ and
$\rho : \R_+ \rightarrow \R_+$ is a bump
function (i.e. a smooth function with $\rho'(0)=0$, $\rho=0$ outside
$[0;1]$ and $\int_{\R_+} \rho(u) du=1$).

\begin{slemma}
\label{propfkl}
The following holds for the mappings $f_{k,l}$~:
\begin{item}{(i)}
It verifies the Lipschitz property
$$  \exists L(k,l)>0 , \ \forall b,b' \in \R^d,
  \forall x,x' \in \overline O_1, \forall z,z' \in \R^{n d_w},
$$

\begin{equation}
\left\vert f_{k,l}(b,x,z) - f_{k,l}(b',x',z') \right\vert \le L(k,l)
\Bigg( (\vert b-b' \vert + \vert x-x' \vert) (1 + \Vert z \Vert +
\Vert z' \Vert)  + \left\Vert z - z' \right\Vert   \Bigg).
\end{equation}
\end{item}

\begin{item}{(ii)}
When $z=0$, it is uniformly bounded as in \eqref{upperboundf}, and
the bound does not depend on $k,l$.
\end{item}
\end{slemma}

\begin{demo}
As soon as $dist(O_1,O) > 1/l_0 \ge 1/l$, $\rho_l(\beta,y,w)=0$ if
$\vert y \vert \ge dist(O_1,O)$; so the integrand vanishes if
$x-y \notin O$ and we can use the properties of $f_k$ on $O$.\\
It is then obvious that $f_{k,l}$ is Lipschitz (in the sense of
\eqref{lipz}) in the variables $b$ and $z$ (and the Lipschitz
constant can be taken the same as the one for $f_k$). Then the only
point to check is for the $x$ variable. But one has
$$ f_{k,l}(b,x,z)
    = \int_{\R^d \times \R^n \times \R^{n d_w}}
        \rho_l((b,x,z)-(\beta,y,w)) f_k(\beta,y,w) d (\beta,y,w);$$
since $x \mapsto \rho_l((b,x,z)-(\beta,y,w))$ is a smooth function
with bounded partial derivatives, we get
\begin{eqnarray*}
\vert D_x f_{k,l}(b,\cdot,z) \vert
&\le&
C_1(k,l) \int_{\Vert (b,x,z)-(\beta,y,w) \Vert \le 1/l}
          \vert f_k(\beta,y,w) \vert d (\beta,y,w) \\
&\le&
C_{k,l} (1+ \Vert z \Vert)
\end{eqnarray*}
where the last inequality results from the linear growth in $z$ of
$f_k$ (using Lemma \ref{lingro}). Thus $x \mapsto f_{k,l}(b,x,z)$
verifies a Lipschitz condition with a constant growing linearly with
$z$, as was to be proved.

Finally, assumption $(ii)$ results from the same property for $f_k$
(($iii$) in Proposition \ref{propfk}) and the Lipschitz property
\eqref{lipz} for $f$~:
\begin{eqnarray*}
\vert f_{k,l}(b,x,0) \vert
&\le&
\int_{\Vert (\beta,y,w) \Vert \le
1/l} \vert f_k(b-\beta,x-y,-w) \vert \rho_l(\beta,y,w) d
(\beta,y,w)\\
&\le& \int_{\Vert (\beta,y,w) \Vert \le 1/l} C_0
 \rho_l(\beta,y,w) d (\beta,y,w) = C_0
\end{eqnarray*}
using the definition of $\rho_l$.
\end{demo}

Now the mappings $f_{k,l}$ are not necessarily pointing outward on
the boundary $\partial \omb$; in order to get this property, we
introduce the new drift
$$g_{k,l}(b,x,z) = f_{k,l}(b,x,z)
     + (\ep_{k,l} + \frac{A}{l}(1+ \Vert z \Vert)) x,$$
where $\ep_{k,l}$ is the constant associated with $\eta=1/l$ in the
uniform continuity assumption $(vi)$ of Proposition \ref{propfk} and
$A$ is a constant which is chosen below.

\begin{sprop}
The mapping $g_{k,l}$ verifies properties $(i)$ and $(ii)$ in Lemma
\ref{propfkl} above. Moreover, there is a positive constant $A$,
depending on $k$ {\bf but not on $l$}, such that $g_{k,l}$ is
pointing outward on the boundary $\partial \omb$.

Therefore, from Theorem 1.4.1 in \cite{blache03} and Theorem 1.3.1
in \cite{blache05}, the BSDE $(M+D)$ with drift $g_{k,l}$ has a
unique solution $(X_{k,l}, Z_{k,l})$ with $X_{k,l} \in \omb$.
\end{sprop}

\begin{demo}
That $g_{k,l}$ verifies $(i)$ and $(ii)$ in Lemma \ref{propfkl}
above is obvious. Now let us prove the last property; consider $x
\in
\partial \omb = S(0,1)$ . If we put $E:=\vert f_{k,l}(b,x,z) -
f_k(b,x,z) \vert$ then
\begin{eqnarray*}
E
&=&
\left\vert \int
        (f_k(b-\beta,x-y,z-w) - f_k(b,x,z))
              \rho_l(\beta,y,w) d(\beta,y,w) \right\vert  \\
&\le&
\int
     \left\vert f_k(b-\beta,x-y,z-w) - f_k(b-\beta,x,z-w) \right\vert
              \rho_l(\beta,y,w) d(\beta,y,w)  \\
&  &
+ \int
     \left\vert f_k(b-\beta,x,z-w) - f_k(b,x,z) \right\vert
              \rho_l(\beta,y,w) d(\beta,y,w)  \\
&\le& \ep_{k,l} + \frac{C}{l} (1+ \Vert z \Vert)
\end{eqnarray*}
where the last inequality is a consequence of $(i)$ and $(vi)$ in
Proposition \ref{propfk} (note that $C$ depends neither on $k$ nor
on $l$).

Since $f_k$ is pointing outward on $\partial \omb$, this shows that
the following lower bound on the radial part of $f_{k,l}$ holds on
$\partial \omb$~:
$$[f_{k,l}(b,x,z)]^{\rm rad} \ge - \ep_{k,l}
- \frac{C}{l} (1+ \Vert z \Vert).$$
Now if we choose $A:=C+1$, which is independent of $k$ and $l$,
then the mapping $g_{k,l}$ is also pointing outward on $\partial
\omb$. This completes the proof of the proposition.
 \end{demo}

\bigskip
{\it Third step.} Now it remains to pass through the limit; first
for $k$ fixed and $l \rightarrow \infty$, then for $k \rightarrow
\infty$. Note that the terminal value $U \in \omb$ is fixed
independently of $k,l$.

\medskip
First let $k$ fixed and $l \rightarrow \infty$. Let $(X^l,Z^l)$ be the
(unique) solution of the equation $(M+D)$ with the drift $g_{k,l}$
defined above. In general, we will drop the superscript $k$ in the
notations, for simplicity. \\
For $l,l' \in \N^*$, we put $\tilde X_t^{l,l'}:=(X_t^l, X_t^{l'})$,
$\tilde Z_t^{l,l'}:=(Z_t^l, Z_t^{l'})$ and $V_t^l:=(B_t^y, X_t^l,
Z_t^l)$. We apply I\^ o's formula between $t$ and $T$ to $\hat \Psi$
defined by
$$
\left\{
\begin{array}{ll}
\hat \Psi = \Psi \approx \de^2
&  \hbox{ if a general connection is used}; \\
\hat \Psi = \de^2
&  \hbox{ in the case of the Levi-Civita connection.}
\end{array}
\right.
$$
We get
\begin{eqnarray}
-\hat \Psi (\tilde X_t^{l,l'})
    &=& \int_t^T D \hat \Psi (\tilde X_s^{l,l'})
            \left( \tilde Z_s^{l,l'} d W_s \right)     \nonumber   \\
    & & + \half \int_t^T  \left(
      \sum_{i=1}^{d_W} {}^t [{}^t \tilde Z_s^{l,l'}]^i
                    \Hess \hat \Psi (\tilde X_s^{l,l'})
                  [{}^t \tilde Z_s^{l,l'}]^i \right) ds    \nonumber \\
    & & + \int_t^T  D \hat \Psi(\tilde X_s^{l,l'})
                                   \left(
                                 \begin{array}{c}
                                     g_{k,l}(V_s^l)  \\
                                     g_{k,l'}(V_s^{l'})
                                 \end{array}
                                   \right) ds
\label{ito}
\end{eqnarray}

In fact this does not include the case of a drift $f$ which does not
depend on $z$; for this simpler case, see the remark at the end of
this section.

 The first term on the right is a martingale; the term involving the
Hessian is bounded below by
$$
\hat \al \int_t^T \left\Vert
  \overset{X_s^{l'}}{\underset{X_s^l}{\Vert}} Z_s^l - Z_s^{l'}
                                         \right\Vert_r^2 ds
      - \hat \beta \int_t^T \hat \Psi(\tilde X_s^{l,l'})
             (\Vert Z_s^l \Vert_r^2 + \Vert Z_s^{l'} \Vert_r^2) ds,
$$
or
\begin{equation}
\al \int_t^T \left\Vert Z_s^l - Z_s^{l'} \right\Vert^2 ds
 - \beta \int_t^T \hat \Psi(\tilde X_s^{l,l'})
             (\Vert Z_s^l \Vert^2 + \Vert Z_s^{l'} \Vert^2) ds,
\label{min1}
\end{equation}
for positive constants $\al$ and $\beta$. The first bound results
from Proposition 2.3.2 in \cite{blache03} if sectional curvatures
are nonpositive (in this case $\beta=0$) and from estimate (4.2) in
\cite{blache05} in the other cases. The second bound is then a
consequence of Proposition \ref{proptp}.\\
The last term involving $g_{k,l}$ and $g_{k,l'}$ can be rewritten as
\begin{eqnarray*}
A_{l,l'} + \int_t^T  D \hat \Psi(\tilde X_s^{l,l'})
                                   \left(
                                 \begin{array}{c}
                                     f_k(V_s^l)  \\
                                     f_k(V_s^{l'})
                                 \end{array}
                                   \right) ds
\end{eqnarray*}
where
\begin{eqnarray}
A_{l,l'}
&:=&
\int_t^T  D \hat \Psi(\tilde X_s^{l,l'})
                                   \left(
                                 \begin{array}{c}
                                     (g_{k,l}-f_k)(V_s^l)  \\
                                     (g_{k,l'}-f_k)(V_s^{l'})
                                 \end{array}
                                   \right) ds.
\label{defAl}
\end{eqnarray} Using $(ii)$ in Proposition
\ref{propfk}, we have
\begin{equation}
\int_t^T  D \hat \Psi(\tilde X_s^{l,l'})
                                   \left(
                                 \begin{array}{c}
                                     f_k(V_s^l)  \\
                                     f_k(V_s^{l'})
                                 \end{array}
                                   \right) ds
\ge - C \cdot k \int_t^T \hat \Psi(\tilde X_s^{l,l'})
              ( 1 + \Vert Z_s^l \Vert + \Vert Z_s^{l'} \Vert) ds
\label{min2}
\end{equation}
for a constant $C$ independent of $l$. Now using \eqref{min1} and
\eqref{min2}, we get by taking the expectation in \eqref{ito}
\begin{eqnarray}
\esp (\hat \Psi (\tilde X_t^{l,l'})) +
\al \esp \int_t^T \left\Vert Z_s^l - Z_s^{l'} \right\Vert^2  ds
&\le&  C^1_k \Bigg( \int_t^T \esp \bigg( \hat \Psi(\tilde X_s^{l,l'})
      (1 + \Vert Z_s^l \Vert_r^2 \nonumber \\
&   & \hphantom{C^1_k \int_t^T \esp \bigg( }
  + \Vert Z_s^{l'} \Vert_r^2) \bigg) ds +\esp \vert A_{l,l'} \vert \Bigg)\nonumber  \\
&\le&  C^2_k \left( \int_t^T \esp \left( \hat \Psi(\tilde
X_s^{l,l'}) \right) ds + \esp \vert A_{l,l'} \vert \right) \nonumber
\\
\label{min3}
 \end{eqnarray}
where the last inequality is obtained using H\"older's inequality,
noting the uniform exponential integrability of the $Z$-processes
(Lemma \ref{Emu}) and that $\hat \Psi^2 \le \hat \Psi$ on the
compact set $\omb \times \omb$. An application of Gronwall's Lemma
gives
\begin{equation}
\esp (\hat \Psi (\tilde X_t^{l,l'})) \le C_k \esp \vert A_{l,l'}
\vert. \label{gronwall}
\end{equation}

\begin{slemma}
The expectation $\esp \vert A_{l,l'} \vert$, where $A_{l,l'}$ is
defined by \eqref{defAl}, tends to zero as $l,l'$ tend to infinity.
\end{slemma}

\begin{demo}
Since $\hat \Psi \approx \delta^2$ and $g_{k,l}(b,x,z) = f_{k,l}(b,x,z)
     + (\ep_{k,l} + \frac{A}{l}(1+ \Vert z \Vert)) x,$ we have
\begin{eqnarray*}
\esp \vert A_{l,l'} \vert
&\le& \int_0^T \esp \left\vert D \hat \Psi
(\tilde X_s^{l,l'})
                                   \left(
                                 \begin{array}{c}
                                     (f_{k,l}-f_k)(V_s^l)  \\
                                     (f_{k,l'}-f_k)(V_s^{l'})
                                 \end{array}
                                   \right) \right\vert ds \\
&&
+ \int_0^T \esp \left\vert D \hat \Psi (\tilde X_s^{l,l'})
                                   \left(
                                 \begin{array}{c}
          (\ep_{k,l} + \frac{A}{l}( 1 + \Vert Z_s^l \Vert)) X_s^l\\
          (\ep_{k,l'} + \frac{A}{l'}( 1 + \Vert Z_s^{l' }\Vert)) X_s^{l'}
                                 \end{array}
                                   \right) \right\vert ds \\              &\le&
\int_0^T  \esp \left\vert D \hat \Psi (\tilde X_s^{l,l'})
                                   \left(
                                 \begin{array}{c}
                                     (f_{k,l}-f_k)(V_s^l)  \\
                                     (f_{k,l'}-f_k)(V_s^{l'})
                                 \end{array}
                                   \right) \right\vert ds \\
&& + C \int_0^T  \esp \left( \delta (\tilde X_s^{l,l'})
  (\ep_{k,l} + \ep_{k,l'} + (\frac{A}{l}+\frac{A}{l'})( 1 + \Vert Z_s^l \Vert + \Vert Z_s^{l'}
    \Vert))  \right) ds;
\end{eqnarray*}
The last inequality is obtained with Lemma \ref{upbddpsi}; the first
term on the right tends to zero because $f_{k,l}$ tends uniformly to
$f_k$ as $l$ tends to infinity (this results from the properties of
convolution and the uniform continuity of $f_k$). The second term
tends to zero because the integral is bounded independently of $l$
(using the uniform exponential integrability condition of Lemma
\ref{Emu}) and $\ep_{k,l}$ tends to zero as $l$ tends to infinity
(see $(vi)$ in Proposition \ref{propfk}).
\end{demo}

As a consequence, as $l,l'$ tend to infinity, $\esp \hat \Psi(\tilde
X_s^{l,l'}) \rightarrow 0$ (in fact we have obviously $\esp \int_0^T
\hat \Psi(\tilde X_s^{l,l'}) ds \rightarrow 0$), which shows that
the sequence of processes $X^l$ is Cauchy in $L^2(\Omega \times
[0;T])$ thus converges to a process $X^k$.
From \eqref{min3}, we deduce that the sequence of processes $Z^l$ is
also Cauchy in  $L^2(\Omega \times [0;T])$ and  thus converges to a
process $Z^k$ as $l$ tends to infinity.\\

At the end, the pair $(X^k, Z^k)$ solves the equation $(M+D)$ with
the drift $f_k$; passing through the limit is an easy adaptation of
the second step in the proof of Proposition 4.1.4 in
\cite{blache03}, using once again the uniform continuity of $f_k$
and the uniform convergence of $f_{k,l}$ to $f_k$. We recall briefly
the proof for the sake of completeness.

For $k$ fixed, let us show that the following expectation tends to
zero as $l$ tends to $\infty$ :
\begin{eqnarray*}
\esp \bigg\vert U & - & \int_t^T Z^k_s d W_s - \int_t^T \left( -
\half
   \Gamma_{ij}(X^k_s) (\left[ Z^k_s \right]^i \vert
         \left[ Z^k_s \right]^j) + f_k(B_s^y,X^k_s, Z^k_s)  \right) d s        \\
  & + & \int_t^T Z^l_s d W_s + \int_t^T \left( - \half
   \Gamma_{ij}(X^l_s) (\left[ Z^l_s \right]^i \vert
         \left[ Z^l_s \right]^j) + f_{k,l}(B_s^y,X_s^l, Z_s^l)  \right) d s
        \bigg\vert.
\end{eqnarray*}
This expectation is bounded above by
\begin{eqnarray*}
&   & \esp\left( \int_0^T \Vert Z_s^l-Z^k_s \Vert^2 ds \right)^\half\\
& + & \esp \left( \int_0^T \vert \Gamma_{ij}(X^k_s) -
\Gamma_{ij}(X_s^l)
     \vert \cdot \vert (\left[ Z^k_s \right]^i \vert
     \left[ Z_s^k \right]^j) \vert ds \right)                    \\
& + & \esp \left( \int_0^T \vert \Gamma_{ij}(X_s^l) \vert
     \left\vert  (\left[ Z_s^k \right]^i \vert \left[ Z_s^k \right]^j)
        -(\left[ Z_s^l \right]^i \vert \left[ Z_s^l \right]^j)
       \right\vert ds \right)                                  \\
& + & \esp \left( \int_0^T \vert f_k(B_s^y,X^k_s, Z^k_s)
                          - f_{k,l}(B_s^y,X_s^l, Z_s^l)  \vert ds \right).
\end{eqnarray*}
We know that the first expectation tends to zero; the second term
tends to zero by dominated convergence (at least for a subsequence
of $(X^l)$, but it doesn't matter since $(X^l)$ is Cauchy). Let
$E_1$ denote the next expectation; then we can write
\begin{eqnarray*}
E_1 & \le & C \esp \left( \int_0^T \Vert Z_s^l - Z^k_s \Vert
    (\Vert Z_s^l \Vert + \Vert Z^k_s \Vert) ds \right)           \\
& \le & \sqrt 2 C \esp \left( \int_0^T \Vert Z_s^l - Z^k_s \Vert ^2
ds
        \right) ^\half
        \esp \left( \int_0^T (\Vert Z_s^l \Vert^2 + \Vert Z_s^k \Vert^2)
             ds \right) ^\half.
\end{eqnarray*}
The first integral tends to zero and the second is bounded because
$Z^l$ converges in $L^2$; hence $E_1$ tends to zero.

Finally, let $E_2$ denote the last integral; then
\begin{eqnarray*}
E_2 &\le& \esp \left( \int_0^T \vert f_k(B_s^y,X^k_s, Z^k_s)
                          - f_{k}(B_s^y,X_s^l, Z_s^l)  \vert ds
                          \right) \\
&  & + \esp \left( \int_0^T \vert f_k(B_s^y,X^l_s, Z^l_s)
                          - f_{k,l}(B_s^y,X_s^l, Z_s^l)  \vert ds
                          \right);
\end{eqnarray*}
now the first term goes to zero as $l \rightarrow \infty$ by
dominated convergence (at least for a subsequence for which we have
a.s. convergence); for the second term, the zero limit results also
from dominated convergence; indeed, it suffices to write
\begin{eqnarray*}
&& \vert f_k(B_s^y,X^l_s, Z^l_s)
- f_{k,l}(B_s^y,X_s^l, Z_s^l)  \vert \le \\
&&\int_{\R^d \times \R^n \times \R^{n d_w}}
        \vert f_k(B_s^y,X^l_s, Z^l_s)-f_k(B_s^y-\beta,X_s^l-y, Z_s^l-w)
         \vert \rho_l(\beta,y,w) d (\beta,y,w)
\end{eqnarray*}
and use Properties $(i)$ and $(vi)$ of Proposition \ref{propfk}.

Hence the limit in $L^1$ of $X_t^l$ is
$$U - \int_t^T Z^k_s d W_s - \int_t^T \left( - \half
   \Gamma_{ij}(X^k_s) (\left[ Z^k_s \right]^i \vert
         \left[ Z_s \right]^j) + f(B_s^y,X^k_s, Z^k_s)  \right) d s;$$
we know that it is also $X^k_t$, so by continuity~:
$$ a.s., \ \forall t, \
X^k_t=U - \int_t^T Z^k_s d W_s - \int_t^T \left( - \half
   \Gamma_{ij}(X^k_s) (\left[ Z^k_s \right]^i \vert
         \left[ Z^k_s \right]^j) + f(B_s^y,X^k_s, Z^k_s)  \right) d s.$$
The proof is complete.

\bigskip
Finally suppose that $k \rightarrow \infty$. In this case we can
follow the proof above, making the changes accordingly; for
instance, the expression \eqref{defAl} becomes now~:
\begin{eqnarray*}
A_{k,k'} &:=& \int_t^T  D \hat \Psi(\tilde X_s^{k,k'})
                                   \left(
                                 \begin{array}{c}
                                     (f_k-f)(V_s^k)  \\
                                     (f_{k'}-f)(V_s^{k'})
                                 \end{array}
                                   \right) ds;
\end{eqnarray*}
remark that \eqref{min3} is replaced by \eqref{monot} so in
particular the constant in the new inequality \eqref{gronwall} does
not depend on $k$. We give the proof of the new lemma we need~:
\begin{slemma}
The expectation $\esp \vert A_{k,k'} \vert$ tends to zero as $k,k'$
tend to infinity.
\label{akk}
\end{slemma}

\begin{demo}
We have
\begin{eqnarray*}
\esp \vert A_{k,k'} \vert &\le&
 C_1 \int_0^T  \esp \left( \delta (\tilde X_s^{k,k'})
  ( \vert f_k(V_s^k)-f(V_s^k) \vert + \vert f_{k'}(V_s^{k'})-f(V_s^{k'}) \vert)
    \right) ds;\\
&\le&
 C \int_0^T  \esp \bigg( \delta (\tilde X_s^{k,k'})
  \big( ( 1+ \Vert Z_s^k \Vert)
  (1_{ \{ \Vert Z_s^k \Vert \ge k \} }+1_{ \{ \vert B_s^y \vert \ge k \}
  })\\
&  & \hphantom{aaaaaaaaaaaaaaaa}
  + (1+\Vert Z_s^{k'} \Vert)
  (1_{ \{ \Vert Z_s^{k'} \Vert \ge k' \} }+1_{ \{ \vert B_s^y \vert \ge k' \}
  }) \big)
 \bigg) ds;
\end{eqnarray*}
The last inequality uses the definition \eqref{deffk} of $f_k$ and
the linear growth of $f$ (Proposition \ref{lingro}). As a
consequence of Markov inequality and the uniform exponential
integrability condition (Lemma \ref{Emu}), one has $1_{ \{ \Vert
Z_s^k \Vert \ge k\}} \rightarrow 0$ when $k \rightarrow \infty$.
Then the convergence to zero follows once again from dominated
convergence.
\end{demo}

As above, this delivers a pair of processes $(X,Z)$. To prove that
this pair solves the equation $(M+D)$, we can yet follow the proof
above, and the only difficulty is to deal with the term involving
$f$. But, putting $V_s:=(B_s^y, X_s, Z_s)$ and $V_s^k:=(B_s^y,
X_s^k, Z_s^k)$ we have
\begin{eqnarray*}
\esp \left( \int_0^T \vert f(V_s) - f_k(V_s^k) \vert ds \right)
&\le&
\esp \left( \int_0^T \vert f(V_s) - f(V_s^k) \vert ds \right) \\
& & + \esp \left( \int_0^T \vert f(V_s^k) - f_k(V_s^k) \vert ds
\right).
\end{eqnarray*}
Then the two terms on the right tend to zero as in the proof of
Lemma \ref{akk}, and this completes the proof when $k \rightarrow
\infty$.

{\bf Remark : } If $f$ depends only on $x$ (and not on $z$), we can
take the convex function $\Psi \approx \delta^p$ of \cite{blache03}.
Then $\Hess \Psi$ is nonnegative and the same reasoning as above
(even much simpler) gives a limit process $(X_k)$ but not $(Z_k)$.
In fact, the only problem to deal with is to verify that $X_k$
solves equation $(M+D)$ with drift $f_k$. But this results from
Section 2.4 in \cite{blache03}, a localization argument and passing
through the limit in a submartingale (using the continuity of $f_k$)
as in the proof of Proposition 4.1.4 in \cite{blache03}.

%%%%%%%%%%%%%%%%%%%%%%%%%%%%%%%%%%%%%%%%%%%%%%%%%%%%%%%%%%%%%%%%%%%%%%%%%
%%%%%%%%%%%%%%%%%%%%%%%%%%%%%%%%%%%%%%%%%%%%%%%%%%%%%%%%%%%%%%%%%%%%%%%%%
%%%%%%%%%%%%%%%%%%%%%%%%%%%%%%%%%%%%%%%%%%%%%%%%%%%%%%%%%%%%%%%%%%%%%%%%%
\section{Random terminal times}
The purpose of this section is to explain how to extend the above
results to the case of a random terminal time. In the Lipschitz
context, this has been done in Section 5.3 in \cite{blache03}.

We consider here the following equation
$$(M+D)_\tau
\left\{
      \begin{array}{l}
          d X_t = Z_t d W_t + \left( - \half \Gamma_{jk}(X_t)
            ([Z_t]^k \vert [Z_t]^j) + f(B_t^y,X_t,Z_t) \right) d t  \\
      X_\tau=U^\tau                                         \\
      \end{array}
    \right.   $$
where $\tau$ is a stopping time with respect to the filtration used
and $U^\tau$ is a $\omb$-valued, ${\cal F}_\tau$-measurable random
variable. It is the counterpart of equation $(M+D)$ on the random
interval $[0; \tau]$.

For bounded random terminal times (i.e. $\tau \le T$ where $T$ is a
deterministic constant), the proofs given in the preceding sections
can be used again, together with Theorem 5.3.1 of \cite{blache03},
which gives existence and uniqueness in the Lipschitz context. Thus
we have the

\begin{stheorem}
\label{existence2} We consider BSDE $(M+D)_\tau$ with $\omb= \{ \chi
\le c \}$ and $\tau \le T$ a.s. Then under the same assumptions and
in the same cases as in Theorem \ref{existence}, this BSDE has a
unique solution $(X,Z)$, with $X \in \omb$.
\end{stheorem}

%%%%%%%%%%%%%%%%%%%%%%% nonbounded stopping time %%%%%%%%%%%%%%%%%%%%%%%%
\bigskip
Next, we consider a stopping time $\tau$ that is finite a.s. and
verifies the exponential integrability condition
\begin{equation}
  \exists \rho > 0 :  \esp (e^{\rho \tau})<\infty.
\label{integst}
\end{equation}
Examples of such stopping times are exit times of uniformly elliptic
diffusions from bounded domains in Euclidean spaces.

As in the Lipschitz case, we need to add restrictions on the drift
$f$; indeed, the proofs above rely heavily on the construction of
the submartingale $(S_t)_t = (\exp(A_t) \Psi(\tilde X_t))_t$, so we
need to keep the integrability of $S_\tau$. Looking at the
computations in the uniqueness part, the conclusion is then the
same~: this integrability holds for "small" drifts, i.e. there is a
constant $h$ with $0<h<\rho$ such that, under the following
condition on the constants in \eqref{lipz}, \eqref{monot} and
\eqref{upperboundf}
\begin{equation}
L < h, \ \ \nu > -h, \ \ L_2 < h,
\label{smalldrift}
\end{equation}
the integrability required holds, so $(S_t)_{0 \le t \le \tau}$ is a
true submartingale.

In fact we can state the following result, whose proof goes exactly
the same as in Section 5.3 in \cite{blache03}.

\begin{stheorem}
\label{existence3} We consider the BSDE $(M+D)_\tau$ with $\tau$ a
stopping time verifying the integrability condition \eqref{integst};
the function $\chi$ used to define the domain $\omb$ is supposed as
usual to be strictly convex. Then under the same assumptions and in
the same cases as in Theorem \ref{existence}, if moreover we suppose
that $f$ is "small" (i.e. verifies condition \eqref{smalldrift}
above), this BSDE has a unique solution $(X,Z)$.
\end{stheorem}

{\bf Remark~:} Of course the applications to PDEs in \cite{blache03}
are still valid under the current assumptions.

\bibliography{biblio}
\bibliographystyle{abbrv}

\end{document}